\documentclass[11pt,fleqn]{article}
\RequirePackage[T1]{fontenc}
\RequirePackage{amsthm,amsmath}
\RequirePackage[square,numbers]{natbib}
\usepackage[english]{babel}
\RequirePackage[colorlinks,citecolor=blue,urlcolor=blue]{hyperref}
\usepackage{enumerate}
\usepackage{amsmath,amsfonts,amssymb,euscript,mathrsfs}
\usepackage{fullpage}
\usepackage{amsmath}
\usepackage{dsfont}
\makeatletter
\def\captionof#1#2{{\def\@captype{#1}#2}}
\makeatother
\def\1{\mbox{\bf 1}}
\def\R{\mathbb{R}}

\def\N{\mathbb{N}}
\def\P{\mathbb{P}}
\def\E{\mathbb{E}}

\def\R{\mathbb{R}}

\def\Z{\mathbb{Z}}

\newtheorem{theo}{Theorem}
\newtheorem{lem}{Lemma}
\newtheorem{prop}{Proposition}

\newtheorem{cor}{Corollary}
\newtheorem{Def/Prop}{Definition-Proposition}

\newcounter{exos}
\renewcommand\theexos{\arabic{exos}}

\newcounter{prob}
\renewcommand\theprob{\arabic{prob}}

\begin{document}
\title{Ergodic properties of some Markov chains models in random environments}
\date{}

\author{
Lionel Truquet \footnote{UMR 9194 CNRS CREST, ENSAI, Campus de Ker-Lann, rue Blaise Pascal, BP 37203, 35172 Bruz cedex, France.}
 }

\maketitle

\begin{abstract}
We study ergodic properties of some Markov chains models in random environments when the random Markov kernels that define the dynamic satisfy some usual drift and small set conditions but with random coefficients. In particular, we adapt a standard coupling scheme used for getting geometric ergodic properties for homogeneous Markov chains to the random environment case and we prove the existence of a process of randomly invariant probability measures for such chains, in the spirit of the approach of \citet{Kifer} for chains satisfying some Doeblin type conditions. We then deduce ergodic properties of such chains when the environment is itself ergodic. Our results complement and sharpen existing ones by providing quite weak and easily checkable assumptions on the random Markov kernels. As a by-product, we obtain a framework for studying some time series   
models with strictly exogenous covariates. We illustrate our results with autoregressive time series with functional coefficients and some threshold autoregressive processes. 
\end{abstract}

\section{Introduction}
Let $E$ and $F$ be two Polish spaces and for any $x\in F$, we consider a Markov kernel $P_x$ on $E$.
We denote by $\mathcal{B}(E)$ the Borel sigma-field of $E$.
Let also $X:=(X_t)_{t\in\Z}$ a stochastic process taking values in $F$
Our aim is to construct stationary Markov chains in a random environment $(Y_t)_{t\in\Z}$ defined by 
\begin{equation}\label{model}
\P\left(Y_t\in A\vert X,Y_{t-1},Y_{t-2},\ldots\right)=P_{X_{t-1}}\left(Y_{t-1},A\right),\quad t\in\Z. 
\end{equation}
The term Markov chain in random environments (MCRE) comes from the fact that the process $(Y_t)_{t\in\N}$ 
defines a time-inhomogeneous Markov chain conditionally on an exogenous process $X$ called the random environment.
Such a definition has already been used in many previous contributions on the topic. See in particular \citet{Cogburn}, \citet{Orey} or \citet{Kifer}.
However, such references are devoted to quite specific structures, with either discrete state spaces or Doeblin type conditions.
Recently, \citet{DNT} investigated the case of observation-driven models with unbounded state spaces. This specific class of models were studied without assuming a small set assumption for the Markov kernels. For some of their examples, their assumptions are not always easy to check or require restrictive conditions on the noise density while the existence of a small set could be more easily obtained.
An interesting and recent contribution to MCRE that satisfy both drift and small set conditions can be found in \citet{Lovas}. 
In particular, the authors proved some law of large numbers for such MCRE when the system is initialized at a given deterministic value. 
The present paper is mainly motivated by time series analysis and we will study another type of problem, the existence of stationary and ergodic solutions for (\ref{model}) which is one of the crucial tricky point to check when studying stationary time series and their statistical inference.
We will also use different type of assumptions for the random Markov kernels, some of them being weaker.

Note that in (\ref{model}), $(X_s)_{s\geq t}$ is automatically independent from $Y_t$ conditionally on \\$\mathcal{F}_{t-1}:=\sigma\left((Y_u,X_u): u\leq t-1\right)$. In econometrics or time series analysis, such a conditional independence assumption is related to a standard notion called strict exogeneity. 
In this context, past information has an influence on the present value $Y_t$ but the latter one will not impact the present or the future values 
of the exogenous process $X$.
See for instance \citet{Chamb} for a discussion about the various equivalences between some exogeneity notions used in econometrics. 

There exist numerous applications to the study of ergodic properties of MCRE. \citet{Lovas} provided examples in queuing theory, Machine learning or for linear autoregressive processes with random coefficients. One can more generally study non-linear autoregressive processes of the form
\begin{equation}\label{appli}
Y_t=f\left(X_{t-1},Y_{t-1},\varepsilon_t\right),
\end{equation} 
with i.i.d. errors $\varepsilon_t$ independent of the process $(X_t)_{t\in\Z}$. One can then extend some classical time series models usually studied using homogeneous Markov chains theory to more realistic models. We believe that such extensions are fundamental since in the applied statistical literature, the use of exogenous covariates for time series analysis is almost systematic. We will discuss applications of our results to some of these models.

Finally, we mention that the conditional distribution given in (\ref{model}) only depends on $X_{t-1}$ and not on $X_t$. 
We choose such a formulation mainly for compatibility with time series analysis, where the present value of the process can be predicted only using past values of the observations. For a theoretical analysis of Markov chains in random environments, one can always replace the process $X$ by the shifted process $(X_{t+1})_{t\in\Z}$ to go back to the situation where the conditional distribution only depend on $X_t$ and not on $X_{t-1}$ so that both formulations are equivalent. Additionally, if the conditional distribution given in (\ref{model}) has to depend on the whole past $X_{t-1},X_{t-2},\ldots$ of the exogenous process, one can simply replace the single variable $X_{t-1}$ by the new random variable $\left(X_{t-1},X_{t-2},\ldots\right)$ and the state space $F$ of the environment by the product space $F^{\N}$ to go back to our initial formulation. 

The paper is organized as follows. In Section \ref{Section2}, we give our main result when the drift and small set conditions are obtained from one iteration of the chain. An extension to more general chains is given in Section \ref{Section3} whereas two examples of autoregressive processes satisfying our conditions are given in Section \ref{Section4}. The proofs of our results can be found in Section \ref{Section5}. Finally an Appendix section \ref{Section6} contains the proof of an important lemma used in the previous section. 

\section{Assumptions and ergodicity result}\label{Section2}
For measurability issues, we impose the two following conditions. For any $A\in\mathcal{B}(E)$, the mapping $(x,y)\mapsto P_x(y,A)$ is measurable. 
Obviously for any pair $(x,y)\in F\times E$, $P_x(y,\dot)$ has to be a probability measure on $E$.
In what follows, we denote by $\mathcal{L}$ the set of measurable mappings $\gamma:F\rightarrow \R_+$ such that $\E\left[\log^{+}(X_0)\right]<\infty$, with $\log^{+}(x)=\log\left(\max(x,1)\right)$ for $x\in\R_+$.  
The following assumptions will be used. 

\begin{description}
\item[A1] The process $(X_t)_{t\in\Z}$ is stationary.
\item[A2] There exist a measurable mapping $V:E\rightarrow (0,\infty)$ such that for all $(x,y)\in E\times F$ and two elements $\lambda$ and $b$ of $\mathcal{L}$ such that 
for all $x\in F$, $P_xV\leq \lambda(x)V+b(x)$. Moreover 
\begin{equation}\label{upperlimit}
\lim\sup_n \prod_{i=1}^n \lambda\left(X_{-i}\right)^{1/n}<1\mbox{ a.s.}
\end{equation}
\item[A3] There exist a measurable mapping $\eta:(0,\infty)\times F\rightarrow (0,1)$ such that for any $R>0$, one can find a probability kernel $\nu_R$ from $F$ to $E$ such that 
$$P_x(y,A)\geq \eta(R,x)\nu_R(x,A),\quad (x,y,A)\in F\times V^{-1}([0,R])\times\mathcal{B}(E).$$
\end{description}

\paragraph{Notes}

\begin{enumerate}
\item
Assumption {\bf A2} imposes a drift condition for the Markov kernels with some varying coefficients $\lambda(\cdot)$ and $b(\cdot)$. The required conditions for these coefficients are quite weak, e.g. existence of logarithmic moments.
Moreover, under Assumption {\bf A1}, condition (\ref{upperlimit}) is automatically satisfied under ergodicity of the environment $X$ and the Lyapunov coefficient condition
$\E\left[\log(\lambda(X_0))\right]<0$. Indeed, in this case, one can use Birkoff's ergodic theorem to get 
$$\prod_{i=1}^n \lambda(X_{-i})^{1/n}=\exp\left(\frac{1}{n}\sum_{i=1}^n \log\left(X_{-i}\right)\right)\rightarrow \exp\left(\E\left[\log(\lambda(X_0))\right]\right)<1.$$
Such a condition is natural since for the simple case of a real-valued autoregressive process with random coefficients 
\begin{equation}\label{comp}
Y_t=a(X_{t-1})Y_{t-1}+\varepsilon_t
\end{equation}
with i.i.d $(X_{t-1},\epsilon_t)$'s, the condition $\E\left[\log(\lambda(X_0))\right]<0$ is known to be optimal for getting existence of a non-anticipative and stationary solution. See \citet{Boug}, Theorem $2.5$. 
\item
In the case of a stationary but not necessarily ergodic environment, one can still apply the ergodic theorem and the condition (\ref{upperlimit}) simply writes as $\E\left[\log(\lambda(X_0))\vert X^{-1}(\mathcal{I})\right]<0$ a.s. where $\mathcal{I}$ is the set of shift-invariant measurable sets $I$ in $E^{\Z}$ such that $\theta^{-1}I=I$, where $\theta:E^{\Z}\rightarrow \E^{\Z}$ denotes the shift operator defined by $\theta x=(x_{t+1})_{t\in\Z}$ and $X^{-1}(\mathcal{I})$ denotes the sigma-algebra $\left\{\{X\in I\}: I\in \mathcal{I}\right\}$. 
It is difficult to find a more explicit condition in the latter case, except if we impose the more restrictive condition $\lambda(X_0)<1$ a.s. which is the classical condition used for studying stability of Markov chains with a deterministic environment. However, the random environment case offers more flexibility by allowing the varying coefficient $\lambda(\cdot)$ to exceed one for some realizations of the exogenous process.  
\item
Condition (\ref{upperlimit}) is generally weaker than the long-time contractivity condition used in \citet{Lovas}. With our notations, the latter condition writes as
\begin{equation}\label{conc}
\overline{\gamma}:=\lim\sup_{n\rightarrow \infty}\E^{1/n}\left(b(X_0)\prod_{t=1}^n\lambda(X_t)\right)<1
\end{equation}
When $b(\cdot)\geq 1$, using stationarity, we deduce that there exist $\kappa\in (0,1)$ such that $\E\left(\prod_{i=1}^n\lambda(X_{-i})\right)\leq \kappa ^n$ when $n$ is large enough.
Using Markov's inequality, we deduce that there exists $\alpha\in (\kappa,1)$ such that $\sum_{n\geq 1}\P\left[\prod_{i=1}^n\lambda(X_{-i})^{1/n}>\alpha\right]<\infty$ and
from the Borel-Cantelli lemma, we get $\lim\sup_n \prod_{i=1}^n\lambda(X_{-i})^{1/n}<1$ a.s. which is equivalent to our condition. 
Note also that when the coordinates of $X$ are independent, (\ref{conc}) reduces to $\E\left[\lambda(X_0)\right]<1$ while our condition writes as $\E\left[\log\left(\lambda(X_0)\right)\right]<0$ which is weaker in general if we use Jensen's inequality.
\item 
Using the classical terminology used for homogeneous Markov chains, Assumption {\bf A3} entails that any level set of the drift function $V^{-1}([0,R])$ is a small set for the kernels $P_x$. This assumption is then more restrictive than the minorization condition $2.4$ of \citet{Lovas} who only assumed existence of small set of the form $V^{-1}\left([0,R(x)]\right)$ where the value $R(x)$ depends on the drift parameters and the constant $\overline{\gamma}$ defined in (\ref{conc}). However, inspection of the proof of Theorem \ref{mainresult} shows that we only require a large value of $R$ for the small set, see (\ref{level1}), but this value depends in a complicated way of the process $X$, this is why we prefer to use {\bf A3}.
On the other hand, we do not impose any specific condition related to the behavior of the minorization function $\eta$ near $0$, such as the assumption $2.6$ used by \citet{Lovas}. As a consequence, for the simple autoregressive process (\ref{comp}), it is straightforward to get {\bf A3} when the noise term $\epsilon_1$ has a positive density lower-bounded on any compact set, while the condition $2.6$ of \citet{Lovas} seems to be difficult to get when $a(\cdot)$ is unbounded and without additional restriction on the noise density.    
\end{enumerate}
We then get the following result.

\begin{theo}\label{mainresult}
Suppose that Assumptions {\bf A1-A3} hold true. There then exists a stationary process $\left((Y_t,X_t)\right)_{t\in\Z}$ satisfying (\ref{model}) and the distribution of such a process is unique. If in addition the process $(X_t)_{t\in\Z}$ is ergodic, the process $\left((Y_t,X_t)\right)_{t\in\Z}$ is also ergodic.
\end{theo}

When the environment is ergodic, our result entails the following strong law of large numbers for the unique stationary solution. If $Z_t=f\left((X_{t+i},Y_{t+i})_{i\in\Z}\right)$ with a measurable mapping $f:(F\times E)^{\Z}\rightarrow \R$ such that $\E\vert Z_0\vert<\infty$, then $\lim_{n\rightarrow \infty}n^{-1}\sum_{t=1}^n Z_t=\E Z_0$ a.s. 

On the other hand, contrarily to \citet{Lovas}, we do not provide a weak law of large numbers when  $Z_t=\phi(Y_t^z)$, where
$Y_t^z$, $t\geq 1$, denotes the iterations of chain (\ref{model}) initialized with $Y_0^z=z$. Such a result could need more 
technical details and we did not investigate it.

One can also show that $\P_{Y_t^z}$ converges in total variation to $\P_{Y_0}$ as $t\rightarrow \infty$ where $\P_{Y_0}$ is the marginal distribution of our stationary solution. We discuss this point just after the statement of Corollary \ref{back}, see (\ref{compar}). However, our assumptions do not help to get a rate of convergence for $\P_{Y_t^z}$ as it is done in \citet{Lovas}. We believe that a precise rate could be more difficult to get only using (\ref{upperlimit}).

\section{Extension to more general chains}\label{Section3}
As in \citet{Lovas}, we now assume that the drift/small set condition is obtained after a given number of iterations of the Markov kernels.
This kind of extension is natural when we face to time-inhomogeneous Markov chains.
In what follows, we recall that the product $T=RS$ of two Markov kernels $R,S$ on $E$ is the Markov kernel defined by $T(y,A)=\int_ER(y,dy')S(y',A)$, $(y,A)\in E\times \mathcal{B}(E)$. 
More precisely, we will assume the existence of a positive integer $p$ such that the following assumptions will be satisfied.

\begin{description}
\item[A4]
There exist measurable mappings $V:E\rightarrow (0,\infty)$ and $\lambda,b:F^p\rightarrow (0,\infty)$ such that for $(x_1,\ldots,x_p)\in F^p$,
$$\left[P_{x_1}\cdots P_{x_p}\right]V\leq \lambda\left(x_p,x_{p-1},\ldots,x_1\right)V+b\left(x_p,x_{p-1},\ldots,x_1\right).$$
Moreover $\E\log^{+}\lambda(X_p,\ldots,X_1)<\infty$, $\E\log^{+}b(X_p,\ldots,X_1)<\infty$ and
\begin{equation}\label{multi}
\lim\sup_{n\rightarrow \infty}\prod_{k=0}^{n-1}\lambda\left(X_{-1-kp},\ldots,X_{-(k+1)^p}\right)^{1/n}<1\mbox{ a.s.}
\end{equation}
\item[A5]
There exist a measurable mapping $\eta:(0,\infty)\times F^p\rightarrow (0,1)$ such that for any $R>0$, one can find a probability kernel $\nu_R$ from $F^p$ to $E$ such that
$$\left[P_{x_1}\cdots P_{x_p}\right](y,A)\geq \eta(R,x_p,\ldots,x_1)\nu(x_p,\ldots,x_1,A),\quad (x,y,A)\in F^p\times V^{-1}([0,R])\times\mathcal{B}(E).$$
\end{description}

We then get the following result which generalizes Theorem \ref{mainresult}.

\begin{theo}\label{mainresult2}
Suppose that Assumptions {\bf A1} and {\bf A4-A5} hold true. There then exists a stationary process $\left((Y_t,X_t)\right)_{t\in\Z}$ satisfying (\ref{model}) and the distribution of such a process is unique. If in addition the process $(X_t)_{t\in\Z}$ is ergodic, the process $\left((Y_t,X_t)\right)_{t\in\Z}$ is also ergodic.
\end{theo}

\paragraph{Note.} Assume $p>1$ and set for $t\in \Z$, $U_t=\left(X_t,\ldots,X_{t-p+1}\right)$. Under the integrability conditions given in {\bf A4} and when the process $X$ is stationary and $p>1$, Birkoff's ergodic theorem implies that 
$$\lim_{n\rightarrow \infty}\prod_{k=0}^{n-1}\lambda\left(U_{-1-kp}\right)^{1/n}=\exp\left(\E\left[\log(U_{-1})\vert X^{-1}\left(\mathcal{I}\right)\right]\right),$$ 
With $\mathcal{I}$ the set of measurable subsets $I$ of $F^{\Z}$ that are invariant for $\theta^p$, i.e. $\theta^{-p}I=I$.
Now, if the process $X$ is only assumed to be ergodic, the conditional expectation appearing in the limit above is not necessarily an expectation.
This differs from the case $p=1$ because a subsequence of an ergodic process is not necessarily ergodic. Let us mention that we obtain the limit $\exp\left(\E\log(U_{-1})\right)$ when 
$X$ is mixing, i.e. for every pair $(A,B)$ of measurable subsets in $F^{\Z}$, we have
\begin{equation}\label{mixing}
\lim_{n\rightarrow \infty}\P\left(X\in A,\theta^n X\in B\right)=\P(X\in A)\P(X\in B).
\end{equation} 
Indeed, if we apply (\ref{mixing}) for $n=kp$ and $A=B=I$ with $I\in\mathcal{I}$, we obtain that $\P(X\in I)$ is necessarily $0$ or $1$.
The mixing property (\ref{mixing}), which is stronger than ergodicity, is satisfied for instance when the process is $\alpha-$mixing.
Hence, when the process is mixing, (\ref{multi}) holds true as soon as $\E\log(U_{-1})<0$. We defer the reader to \citet{sam}, Chapter $2$, for a nice and concise introduction about the ergodicity and mixing concepts and with a summary of the various links existing between these notions.

\section{Examples of autorgressive processes}\label{Section4}

In this section, we present two classical examples of autoregressive processes for which the stability properties are usually 
established using Markov chains techniques. We directly present a version with exogenous regressors and show that our results can be used to 
study these extensions. Let us mention that our examples are mainly illustrative and we argue that our results allow to extend most of the non-linear autoregressive time series models usually studied with Markov chain techniques. We defer the reader to \citet{DoucMoulines}, Section $6.2.4$, for additional examples. 
In these two examples, we consider two independent stochastic processes $(X_t)_{t\in\Z}$ and $\left(\varepsilon_t\right)_{t\in\Z}$ taking values respectively in $F:=\R^d$ and $E:=\R$. We further assume the process $(X_t)_{t\in\Z}$ is stationary and ergodic and that the $\varepsilon_t'$s are i.i.d. The class of function $\mathcal{L}$ was defined at the beginning of Section \ref{Section2}.

\subsection{Threshold autoregressive process}
Let $a_i,b_i,r:\R^d\rightarrow \R$, $i=1,2$, five measurable functions. Set for $t\in\Z$,
\begin{equation}\label{ex1}
Y_t=\left(b_1(X_{t-1})+a_1(X_{t-1})Y_{t-1}\right)\mathds{1}_{Y_{t-1}\leq r(X_{t-1})}+\left(b_2(X_{t-1})+a_2(X_{t-1})Y_{t-1}\right)\mathds{1}_{Y_{t-1}> r(X_{t-1})}+\varepsilon_t.
\end{equation}
Set $\lambda(x)=\max\left(\vert a_1(x)\vert,\vert a_2(x)\vert\right)$.

\begin{prop}\label{example1}
Assume that $\E\vert\varepsilon_0\vert<\infty$, the distribution $\varepsilon_0$ has a positive density $f$ lower-bounded on any compact subset of $\R$, $a_i,b_i\in \mathcal{L}$ for $i=1,2$ and that $\E\log\lambda(X_0)<0$. There then exists a unique stationary and ergodic solution $\left((Y_t,X_t)\right)_{t\in\Z}$ to equations (\ref{ex1}). 
\end{prop}

\paragraph{Note.} When $a_i,b_i,r$ are deterministic, a process defined by (\ref{ex1}) is well known and called threshold autoregressive process.
See \citet{tong} or \citet{Tsay}. Extension to a modeling with exogenous covariate has recently been investigated by \citet{DNT}, see their Section $4.4$.
However, we use here a much weaker assumption on the noise density and then substantially improve their stationarity result for this model.

\subsection{Functional coefficients autoregressive time series}

We consider the following model which is an extension of the functional coefficients autoregressive model of \citet{chen}, see also \citet{Fan}. Let $p$ be a positive integer.
\begin{equation}\label{ex2}
Y_t=\sum_{j=1}^pa_j\left(X_{t-1},Y_{t-1},\ldots,Y_{t-p}\right)Y_{t-j}+\varepsilon_t,\quad t\in\Z.
\end{equation}
Let us first introduce some notations. 
For $j=1,\ldots,p$ and $x\in\R^d$, set $b_j(x):=\sup_{y_1,\ldots,y_p\in\R}\vert a_j(y_1,\ldots,y_p,x)\vert$ which is assumed to be finite and
$$A(x)=\begin{pmatrix}b_1(x)&\cdots &b_p(x)\\ \multicolumn{2}{c}{I_{p-1}}& 0_{p-1,1}\end{pmatrix},$$
with $I_{p-1}$ being the identity matrix of size $p-1$ and $0_{p-1,1}$ the null column vector of size $p-1$.  
Note that $A(x)$ is the companion matrix associated to polynomial $\mathcal{P}_x(L)=1-\sum_{j=1}^p b_j(x)L^j$.
Finally let
$$\gamma(X)=\inf_{n\geq 1}\frac{1}{n}\E\left(\log\Vert A(X_{-1})\cdots A(X_{-n})\Vert\right),$$
where $\Vert\cdot\Vert$ denotes an arbitrary norm on the space of real matrices of dimension $p\times p$. 

\begin{prop}\label{example2}
Assume that $\E\vert\varepsilon_0\vert<\infty$, the distribution $\varepsilon_0$ has a positive density $f$ lower-bounded on any compact subset of $\R$
and that $\gamma(X)<0$. Assume furthermore that the process $X$ is mixing in the sense of (\ref{mixing}). There then exists a unique stationary and ergodic solution $\left((Y_t,X_t)\right)_{t\in\Z}$ to equations (\ref{ex1}). 
\end{prop}

\paragraph{Note.} $\gamma(X)$ is called Lyapunov exponent of the sequence of random matrices $\left(A(X_t)\right)_{t\in\Z}$. It is not straightforward to get a more explicit condition for the negativity of this coefficient. For a general sequence of random matrices $\left(A(X_t)\right)_{t\in\Z}$, negativity of the Lyapunov exponent is a classical condition used for defining stationary solution of random affine transformations on $\R^p$, $Z_t=A(X_{t-1})Z_{t-1}+b_t$ where $\left((X_t,b_t)\right)_{t\in\Z}$ is a stationary process. See \cite{Boug}. A more explicit sufficient condition can be obtained if $c_j:=\sup_{x\in\R^d}b_j(x)$ satisfies $\sum_{j=1}^p c_j<1$, since in this case the spectral radius of the companion matrix associated to $\mathcal{P}(L)=1-\sum_{j=1}^p c_j L^j$ is less than one. Note that $\sum_{j=1}^p c_j<1$ is the classical condition used in the model without covariates. See \citet{DoucMoulines}, p. $186$. 

\section{Proof of the results}\label{Section5}

Our approach for proving Theorem \ref{mainresult} is inspired by that of \citet{Kifer}. Denoting $\xi_t=(X_{t-j})_{j\geq 0}$, we prove the existence of some random probability measures $\left(\pi_{\xi_t}\right)_{t\in\Z}$ such that for any $t\in\Z$,
\begin{equation}\label{rendinv}
\pi_{\xi_{t-1}}P_{X_t}:=\int \pi_{\xi_{t-1}}(dy)P_{X_t}(y,\cdot)=\pi_{\xi_t}\mbox{ a.s.}`
\end{equation}

To define such random probability measures, we impose that the mapping $(p,A)\mapsto \pi_p(A)$ satisfies the definition of a probability kernel from $F^{\N}$ to $E$. In \citet{Stenflo}, random probability measures satisfying (\ref{rendinv}) are called randomly invariant.  

Following \citet{Kifer}, natural candidates for $\pi_{\xi_t}$ are given by the almost sure limits of the backward iterations $\delta_yP_{X_{t-n}}\cdots P_{X_t}$ when $n\rightarrow \infty$ where $\delta_y$ denotes the Dirac mass at point $y$. For simplicity of notations, set $Q_{t,n}^{\omega}=P_{X_{t-n}(\omega)}\cdots P_{X_t(\omega)}$.
We remind that the product $RS$ of two Markov kernels $R$ and $S$ on $E$ is the Markov kernel on $E$ defined by $RS(y,A)=\int R(y,dz)S(z,A)$ for $y\in E$ and $A\in\mathcal{B}(E)$. 
To prove existence of such a limit, a possible approach is to get a control of the total variation distance
\begin{equation}\label{sup}
\sup_{A\in \mathcal{B}(E)}\left\vert \delta_zQ_{0,n}^{\omega}(A)-\delta_{z'}Q_{0,n}^{\omega}(A)\right\vert
\end{equation}
for two initial state values $z$ and $z'$ in $E$ and a fixed $\omega\in \Omega$. Note that by stationarity of the process $X$, it is only necessary to prove 
the existence of $\pi_{\xi_0}$.  
Such a path-by-path control will obtained from a coupling argument detailed below. In the rest of the section, we assume that Assumptions {\bf A1-A3} are satisfied.

\subsection{Coupling strategy}

For a given $\omega\in\Omega$ and a positive integer $n$, we define a probability measure $\overline{P}_{n,\omega}$ on $G=(E\times E)^{\{-n,-n+1,\ldots \}}$ in the following way. 
First, for $t\geq -n$, we denote by $Y_t$ and $\overline{Y}_t$ the coordinate mappings, i.e.
$$Y_t\left((y_{-n+j},\overline{y}_{-n+j})_{j\geq 0}\right)=y_t,\quad \overline{Y}_t\left((y_{-n+j},\overline{y}_{-n+j})_{j\geq 0}\right)=\overline{y}_t.$$
We then assume that $\overline{P}_{n,\omega}\left(Y_{-n}=z,\overline{Y}_{-n}=\overline{z}\right)=1$ and for $t\geq -n+1$, we define this probability measure as the distribution of a Markov chain defined as follows. For two real numbers $a$ and $b$, we set $a\vee b=\max(a,b)$. In what follows, we consider a positive real number $R=R_{\omega}$ that will be chosen latter (see the formula (\ref{level1}) given below). 
\begin{itemize}
\item
On the event $\{Y_{t-1}=\overline{Y}_{t-1}\}$, we set 
$$\overline{P}_{n,\omega}\left(Y_t\in A,\overline{Y}_t\in \overline{A}\vert Y_{t-1},\overline{Y}_{t-1}\right)=P_{X_{t-1}(\omega)}(Y_{t-1},A\cap\overline{A}\cdot).$$
\item
On the event $\left\{Y_{t-1}\neq \overline{Y}_{t-1}, V(Y_{t-1})\vee V(\overline{Y}_{t-1})>R\right\}$, we set 
$$\overline{P}_{n,\omega}\left(Y_t\in A,\overline{Y}_t\in \overline{A}\vert Y_{t-1},\overline{Y}_{t-1}\right)=P_{X_{t-1}(\omega)}\left(Y_{t-1},A\right)P_{X_{t-1}(\omega)}\left(Y_{t-1},A\right).$$
\item
Finally, on the event $\left\{Y_{t-1}\neq \overline{Y}_{t-1}, V(Y_{t-1})\vee V(\overline{Y}_{t-1})\leq R\right\}$, we set
\begin{eqnarray*} 
\overline{P}_{n,\omega}\left(Y_t\in A,\overline{Y}_t\in \overline{A}\vert Y_{t-1},\overline{Y}_{t-1}\right)&=&\eta(R,X_{t-1}(\omega))\nu_R(X_{t-1}(\omega),A\cap \overline{A})\\
&+&\left(1-\eta(R,X_{t-1}(\omega))\right)Q_{X_{t-1}(\omega)}(Y_{t-1},A)Q_{X_{t-1}}(\overline{Y}_{t-1},\overline{A}),
\end{eqnarray*}
\end{itemize}
where for $x\in F$, $A\in \mathcal{B}(E)$ and $y\in E$ such that $V(y)\leq R$,
$$Q_x(y,A)=\frac{P_x(y,A)-\eta(R,x)\nu_R(x,A)}{1-\eta(R,x)}.$$
This coupling scheme is classical for getting some bounds for geometric ergodicity of homogeneous Markov chains. See for instance \citet{DM} or \citet{Ro} who applied such a technique to Markov chains satisfying both a drift and a small set condition. 
Let us mention that some bounds for controlling some quantities similar to (\ref{sup}) are also available for time-inhomogeneous Markov chains. However, one cannot use them in our context because the drift parameters considered in \citet{DM} are assumed to be less than one, which is not necessarily the case for the varying parameter $\lambda(\cdot)$ we use in the present paper.

Let us give an interpretation of the proposed coupling scheme and for simplicity. First note that $\overline{P}_{n,\omega}$ is the distribution of a non-homogeneous Markov chain such that under this probability measure the two coordinate processes $(Y_t)_{t\geq -n}$ and $(\overline{Y}_t)_{t\geq -n}$
are both time-inhomogeneous Markov chains with successive transition kernels $P_{X_{-n}(\omega)},P_{X_{-n+1}(\omega)},\ldots$. 
In the first point, we see that when $Y_{t-1}$ equals to $\overline{Y}_{t-1}$, the two next states also coincide and the common next state is simulated with the Markov kernel $P_{X_{t-1}(\omega)}$. 
When the two previous states are different, two situations can occur. As explained in the second point, when $Y_{t-1}$ or $\overline{Y}_{t-1}$ is outside the small set, the two next states are simulated independently from each other with the same transition kernel $P_{X_{t-1}(\omega)}$. On the other hand, when both previous states are inside the small set, the two next states are equal with probability $\eta(R,X_{t-1}(\omega))$ and the common next state is a realization of the dominating measure $\nu_R\left(X_{t-1}(\omega),\cdot\right)$ or, with probability $1-\eta(X_{t-1}(\omega))$, the two next states are simulated independently form each other and with the same Markov kernel $R_{X_{t-1}(\omega)}$. 

Going back to our original goal, one can note that 
$$\sup_{A\in \mathcal{B}(E)}\left\vert \delta_zQ_{0,n}^{\omega}(A)-\delta_{z'}Q_{0,n}^{\omega}(A)\right\vert\leq P_{n,\omega}\left(Y_0\neq \overline{Y}_0\right)$$
due to dual expression of the total variation distance in term of coupling, i.e. for two probability measures $\mu$ and $\mu'$ on $E$,
$$\sup_{A\in\mathcal{B}(E)}\left\vert \mu(A)-\mu'(A)\right\vert=\inf\left\{\P\left(U\neq U'\right): U\sim \mu, U'\sim \mu'\right\}.$$

\subsection{Proof strategy}\label{scheme}

In the time-homogeneous case, the previous coupling approach can be used in the following way. See in particular \citet{Ro} for a more detailed discussion and specific results.
We then first assume that the coefficients $\eta,\lambda$ and $b$ are deterministic and we write $\overline{P}_n$ instead of $\overline{P}_{n,\omega}$. Let $T_i$, $i\geq 1$, the successive random times (starting here from time $t=-n$) such that $V\left(Y_{T_i}\right)\vee V\left(\overline{Y}_{T_i}\right)\leq R$. 
We have for an arbitrary integer $m <n$, 
$$\overline{P}_n\left(Y_0\neq \overline{Y}_0\right)\leq \overline{P}_n\left(T_m\geq n\right)+\P\left(T_m<n, Y_0\neq \overline{Y}_0\right).$$
Since on the event $\{T_m\leq n\}$, we have a probability greater than $\eta^m$ to get a coalescence of the path, we deduce that 
$$\overline{P}_n\left(Y_0\neq \overline{Y}_0\right)\leq \overline{P}_n\left(T_m\geq n\right)+\eta^m.$$
It then remains to bound the probability $\overline{P}_n\left(T_m\geq n\right)$ which can be obtained from the drift condition if $R$ is large enough.
However, in the case of random environments, there are substantial difficulties due to the functions $\lambda(\cdot),b(\cdot),\eta(\cdot)$ that take either very large or very small values for some time points, depending on the environment. This is why we use an approach comparable to that of \citet{DNT}. In particular, we will consider some random time points, only depending on the environment and for which the function $\eta$ remains lower bounded. Moreover,
these successive random time points are sufficiently spaced, so that the drift parameters of the corresponding subsampled chain remain under control.
The effect of the coupling will be then analyzed only at these random time points. The goal of the next subsection is to introduce such random time points.   

\subsection{Control of the random environment}

Our aim here is to define suitable random times only depending on the process $X$. The following result will be central for this goal.
We denote by $\N^{*}$ the set of positive integers. In what follows, we denote by $\overline{E}_{n,\omega}$ the mathematical expectation associated to $\overline{P}_{n,\omega}$.

\begin{prop}\label{crucial}
There exist two random variables $C_1,C_2:\Omega\rightarrow \N^{*}$ and an increasing sequence of random times $(\tau_i)_{i\in\Z}$, $\tau_i:\Omega\rightarrow \Z$ such that the following statements are valid.
\begin{enumerate}
\item $\tau_{-1}\leq -1$, $\tau_0\geq 0$ and for $i\in\Z$, $\tau_i-\tau_{i-1}\geq C_1$, $\P-$a.s. 
\item 
If $\omega\in\Omega$, let $i\in\Z$ and $s\in \N^{*}$ such that $s\geq C_1(\omega)$. 
We then have
$$\overline{E}_{n,\omega}\left[V\left(Y_{\tau_i(\omega)}\right)\vert Y_{\tau_i(\omega)-s}\right]\leq\left(1-1/C_1(\omega)\right)V\left(Y_{\tau_i(\omega)-s}\right)+C_1(\omega),$$
$$\overline{E}_{n,\omega}\left[V\left(\overline{Y}_{\tau_i(\omega)}\right)\vert \overline{Y}_{\tau_{i-1}(\omega)-s}\right]\leq\left(1-1/C_1(\omega)\right)V\left(\overline{Y}_{\tau_{i-1}(\omega)-s}\right)+C_1(\omega).$$
\item
Setting $R=4C_1(C_1+1)$, we have $\eta\left(R,X_{\tau_i}\right)\geq 1/(C_2+1)$, $\P-$a.s.
\item
$\lim_{i\rightarrow \infty}\tau_i=\infty$ and $\lim_{i\rightarrow -\infty}\tau_i=-\infty$. Moreover if $L_n=\sup\left\{i\geq 1: \tau_{-i}\geq -n\right\}$, then
$$\lim_{n\rightarrow \infty}\frac{L_n}{n}>0\quad \P\mbox{-a.s.}$$ 
\end{enumerate}
\end{prop}

To prove Propostion \ref{crucial}, we first state a lemma.
Let $C=(C_1,C_2)$ be a pair of positive integers and $A_C=A_{1,C_1}\cap A_{2,C}$ with 
$$A_{2,C}=\left\{x\in F^{\N}: \eta(2C_1(2C_1+1),x_0)\geq 1/(C_2+1)\right\}$$
and 
$$A_{1,C_1}=\left\{x\in F^{\N}: \sup_{j\geq C_1}\prod_{i=1}^j\lambda(x_i)\leq 1-1/C_1, b(x_1)+\sum_{i\geq 2}\prod_{k=1}^{i-1}\lambda(x_k)b(x_i)\leq C_1\right\}.$$
Clearly, $A_C$ is an element of the sigma-algebra generated by the cylinders set on $F^{\N}$ and from Birkoff's ergodic theorem, we have 
$$\lim_{n\rightarrow \infty}\frac{1}{n}\sum_{i=0}^n\mathds{1}_{A_C}(\xi_i)=\lim_{n\rightarrow \infty}\frac{1}{n}\sum_{i=1}^n \mathds{1}_{A_C}(\xi_{-i})=\P\left(\xi_0\in A_C\vert X^{-1}(\mathcal{I})\right)\mbox{$\P$-a.s.}$$

\begin{lem}\label{lem1}
The following assertions hold true.
\begin{enumerate}
\item
$\P\left(\xi_0\in \cup_{C_1\geq 1}A_{1,C_1}\right)=1$.
\item
Set $\rho_C=\P\left(\xi_0\in A_C\vert X^{-1}(\mathcal{I})\right)$. There exists a pair of positive, integer-valued random variables $C=\left(C_1,C_2\right)$ such that for $\P-$almost all $\omega\in\Omega$, 
$$\lim_{n\rightarrow \infty}\frac{1}{n}\sum_{i=0}^n\mathds{1}_{A_{C(\omega)}}(\xi_i(\omega))=\lim_{n\rightarrow \infty}\frac{1}{n}\sum_{i=1}^n \mathds{1}_{A_{C(\omega)}}(\xi_{-i}(\omega))=\rho_{C(\omega)}(\omega)>0.$$
\end{enumerate}
\end{lem}

\paragraph{Proof of Lemma \ref{lem1}}
\begin{enumerate}
\item
Let us prove the first point.
Let $\omega\in \Omega$. Condition (\ref{upperlimit}) guarantees the existence, for $\P-$almost $\omega\in\Omega$ of a positive integer $\widetilde{C}$ such that 
\begin{equation}\label{stable}
\sup_{j\geq \widetilde{C}}\lambda(X_{-1}(\omega))\cdots\lambda(X_{-j}(\omega))\leq 1-1/\widetilde{C}.
\end{equation}
Let $C_3$ be the first positive integer such that (\ref{stable}) occurs.
Next, we show that for $\P-$almost $\omega$, there exists a positive integer $C_4$ such that
\begin{equation}\label{finite}
b(X_{-1}(\omega))+\sum_{i\geq 2}\lambda(X_{-1}(\omega))\cdots \lambda(X_{-i+1}(\omega))b(X_{-i}(\omega))\leq C_4.
\end{equation}
To show (\ref{finite}), we use the Cauchy criterion. First, from the log-moment assumption on $b(X_0)$ and the Borel-Cantelli lemma, we note that $\lim_{i\rightarrow \infty}b(X_{-i})^{1/i}=1$ a.s. and using (\ref{upperlimit}), we also get 
$$\lim\sup_{i\rightarrow \infty}\left[\lambda(X_{-1})\cdots \lambda(X_{-i+1})b(X_{-i})\right]^{1/i}<1\mbox{ a.s.}$$
which yields to (\ref{finite}). By taking $C_1=\max(C_3,C_4)$, we see that for $\P-$almost $\omega\in\Omega$, there exists a positive integer $C_1$
such that $\xi_0(\omega)\in A_{1,C_1}$ and the first point of the lemma follows. 
\item
Since the sequence of sets $A_C$ is increasing with respect to $C_2$ and a fixed $C_1$ and the function $\eta$ is positive, we have a.s.
\begin{eqnarray*}
\lim_{C_1\rightarrow \infty}\lim_{C_2\rightarrow\infty}\rho_C&=&\lim_{C_1\rightarrow \infty}\P\left(\xi_0\in A_{1,C_1}\vert X^{-1}(\mathcal{I})\right)\\
&=& \P\left(\xi_0\in \cup_{C_1\geq 1}A_{1,C_1}\vert X^{-1}(\mathcal{I})\right)=1,
\end{eqnarray*}
where the last equality follows from the first point.
Hence, for $\P-$almost every $\omega$, there exists a pair $C(\omega)$ of positive integers such that $\rho_{C(\omega)}(\omega)>0$.
Indeed, if it was not the case, there would exist a subset $\overline{\Omega}$ of $\Omega$ with positive probability and such that for any $\omega\in\overline{\Omega}$ and any pair of positive integers $C$, $\rho_C(\omega)=0$, which contradicts the limiting property given just above.
One can always select $C(\omega)$ in such a way that $\omega\mapsto C(\omega)$ defines a random variable. To this end, one can simply take an ordering 
of $\N^{*}\times \N^{*}$ and that the first pair of integers $C(\omega)$ such that $\rho_{C(\omega)}(\omega)>0$.
Moreover, since the set of pair of positive integers is countable, there exists an event $\widetilde{\Omega}$ with probability $1$ such that for any pair of positive integers $C$ and $\omega\in \widetilde{\Omega}$, 
$$\lim_{n\rightarrow \infty}\frac{1}{n}\sum_{i=0}^n \mathds{1}_{A_C}\left(\xi_i(\omega)\right)=\lim_{n\rightarrow \infty}\frac{1}{n}\sum_{i=1}^n \mathds{1}_{A_C}\left(\xi_{-i}(\omega)\right)=\rho_C(\omega).$$
Such a limit is then also valid for $C=C(\omega)$ which proves the result.$\square$
\end{enumerate}

We now proceed to the proof of Proposition \ref{crucial}.

\paragraph{Proof of Proposition \ref{crucial}}

We define the successive random times $0\leq \widetilde{\tau}_0<\widetilde{\tau}_1<\cdots$ and $-1\geq \widetilde{\tau}_{-1}>\widetilde{\tau}_{-2}>\cdots$ such that $\xi_{\widetilde{\tau}_i(\omega)}\in A_{C(\omega)}$ for $\P-$almost $\omega\in\Omega$, with $\omega \mapsto C(\omega)$ being the random variable defined in the statement of Lemma \ref{lem1}. 
We next define the sequence of random times $\left(\tau_i\right)_{i\in\Z}$ as follows. We set $\tau_{-i}=\widetilde{\tau}_{1-(i-1)C_1}$ for any positive integer $i$ and $\tau_i=\widetilde{\tau}_{1+(i+1)C_1}$ for any non-negative integer $i$.
Note that we have $\tau_i-\tau_{i-1}\geq C_1$ a.s. for any integer $i$.
Setting  
$$M_n=\max\left\{i\geq 1: \widetilde{\tau}_{-i}\geq -n\right\},$$
we note that $M_n(\omega)=\sum_{i=1}^n \mathds{1}_{A_{C(\omega)}}\left(\xi_{-i}(\omega)\right)$ and from the second point of Lemma \ref{lem1}, we get 
\begin{equation}\label{exponent}
\lim\inf_{n\rightarrow \infty}\frac{M_n}{n}>0\mbox{ a.s.}
\end{equation}
Now set 
$$L_n(\omega)=\max\left\{i\geq 1: \tau_{-i}(\omega)\geq -n\right\}.$$
Note that $L_n(\omega)$ is simply the integer part of $\frac{M_n(\omega)-1}{C_1(\omega)}$ and we also have
\begin{equation}\label{exponent2}
\lim\inf_{n\rightarrow \infty}\frac{L_n}{n}>0\mbox{ a.s.}
\end{equation} 
We then get the points $1.$, $3.$ and $4.$ of Proposition \ref{crucial}. It remains to show the second point.
For $\omega\in\Omega$ and $s\in\N^{*}$ such that $s\geq C_1(\omega)$, we have 
\begin{eqnarray*}
\overline{E}_{n,\omega}\left[V\left(Y_{\tau_i(\omega)}\right)\vert Y_{\tau_i(\omega)-s}\right]
 &\leq& \prod_{j=1}^s\lambda\left(X_{\tau_i(\omega)-j}(\omega)\right)V\left(Y_{\tau_i(\omega)-s}\right)+
 b\left(X_{\tau_i(\omega)-1}(\omega)\right)\\
 &+&\sum_{k\geq 2}\prod_{j=1}^{k-1}\lambda\left(X_{\tau_i(\omega)-j}(\omega)\right)b\left(X_{\tau_i(\omega)-k}(\omega)\right)\\
 &\leq& \left(1-1/C_1(\omega)\right)V\left(Y_{\tau_i(\omega)-s}\right)+C_1(\omega),
\end{eqnarray*}
where we used the definition of the random variables $\tau_i$ and $C_1$ and the inequality $s\geq C_1$ for getting the last inequality.
The proof for $\overline{Y}$ is identical and the proof of the proposition is now complete.$\square$

\subsection{Convergence of the backward iterations of the chain}
We now fix the value $R$ for the small set of the chain. For $\omega\in\Omega$, let
\begin{equation}\label{level1}
R(\omega)=2C_1(\omega)\left(2C_1(\omega)+1\right).
\end{equation}

To evaluate the effect of the coupling, we now define two Markov chains $\left(Z^{n,\omega}_i\right)_{i\geq 0}$ and $\left(Z^{n,\omega}_i\right)_{i\geq 0}$
as follows. We set $Z^{n,\omega}_0=z$, $\overline{Z}^{n,\omega}_i=\overline{z}$ and for $i\geq 1$, 
$$Z^{n,\omega}_i=Y_{\tau_{-L_n(\omega)+i}},\quad \overline{Z}^{n,\omega}_i=\overline{Y}_{\tau_{-L_n(\omega)+i}}.$$
From Proposition \ref{crucial}, the two Markov chains introduced above will satisfy a drift condition with non-time varying coefficients that remain in a predetermined stability region. Moreover, note that the probability to get a coalescence of the two paths at any time of the form $\tau_i(\omega)+1$, $i\in\Z$, is lower bounded by a random constant $1/\left(C_2(\omega)+1\right)$ not depending on the integer $i$. 
We also set $W^{n,\omega}_i=V\left(Z^{n,\omega}_i\right)+V\left(\overline{Z}^{n,\omega}_i\right)$. 
We summarize the important properties of these two Markov chains below. 

\begin{lem}\label{lem2}
For $\P-$almost all $\omega\in \Omega$ and $i\geq 1$, we have 
$$\overline{E}_{n,\omega}\left[V\left(Z^{n,\omega}_i\right)\vert Z^{n,\omega}_{i-1}\right]\leq\left(1-1/C_1(\omega)\right)V\left(Z^{n,\omega}_{i-1}\right)+C_1(\omega),$$
$$\overline{E}_{n,\omega}\left[V\left(\overline{Z}^{n,\omega}_i\right)\vert \overline{Z}^{n,\omega}_{i-1}\right]\leq\left(1-1/C_1(\omega)\right)V\left(\overline{Z}^{n,\omega}_{i-1}\right)+C_1(\omega)$$
and 
$$\overline{E}_{n,\omega}\left[W^{n,\omega}_i\vert W^{n,\omega}_{i-1}\right]\leq\left(1-1/C_1(\omega)\right)W^{n,\omega}_{i-1}+2C_1(\omega).$$
\end{lem}

\paragraph{Proof of Lemma \ref{lem2}}
The two first assertions of the lemma are a consequence of Proposition \ref{crucial} and the last one a straightforward consequence of the two first ones.$\square$

Next, we will examine the successive return times of the chain $\left(Z^{n,\omega}_i,\overline{Z}^{n,\omega}_i\right)$ inside the set $\left\{(y,\overline{y}): V(y)+V(\overline{y})\leq R(\omega)\right\}$. Note that at such return times, both chains go back to the small set $\left\{V\leq R(\omega)\right\}$ and the probability to fasten the two paths just after such times is lower bounded by $1/C_2(\omega)$. 
More precisely, let 
$$\rho_1^{n,\omega}=\inf\left\{i\geq 0: V\left(Z^{n,\omega}_i\right)+V\left(\overline{Z}^{n,\omega}_i\right)\leq R(\omega)\right\}$$
and for an integer $j\geq 2$, 
$$\rho_j^{n,\omega}=\inf\left\{i>\rho^{n,\omega}_{i-1}: V\left(Z^{n,\omega}_i\right)+V\left(\overline{Z}^{n,\omega}_i\right)\leq R(\omega)\right\}.$$

Set for $i\geq 0$, $\mathcal{H}^{n,\omega}_i=\sigma\left((Z^{n,\omega}_j,\overline{Z}^{n,\omega}_j): 0\leq j\leq i\right\}$.
We now use a classical result for bounding some exponential moments of these return times. See for instance \citet{DN}, Lemma $3.1$. For completeness  
of this work and since their proof only consider time-homogeneous Markov chains, we give a proof in the Appendix section using our notations.

\begin{lem}\label{lem3}
Set $\eta=\frac{2}{2-1/C_1}$ and $D=1+(1-1/C_1)R+2C_1$. We have the two following bounds.
\begin{enumerate}
\item
If $V(z)+V(\overline{z})>R(\omega)$, we have $E_{n,\omega}\left(\eta(\omega)^{\rho_1^{n,\omega}}\right)\leq V(z)+V(\overline{z})$.
\item
For any $j\geq 0$,
$$\overline{E}_{n,\omega}\left(\eta(\omega)^{\rho^{n,\omega}_{j+1}-\rho^{n,\omega}_j}\vert \mathcal{H}^{n,\omega}_{\rho^{n,\omega}_j}\right)\leq D(\omega)\eta(\omega).$$ 
\end{enumerate}
\end{lem}

From Lemma \ref{lem3}, we now deduce the following important result.

\begin{prop}\label{prop1}
There exist two random variables $F$ and $\kappa$ taking values in $(0,\infty)$ and $(0,1)$ respectively and such that
for $\P-$almost all $\omega\in\Omega$, 
$$\overline{P}_{n,\omega}\left(Y_0\neq \overline{Y}_0\right)\leq F(\omega)\left(1+V(z)+V(\overline{z})\right)\kappa(\omega)^n.$$
\end{prop}

\paragraph{Proof of Proposition \ref{prop1}}
Proceeding as in Subsection \ref{scheme}, we have the bound
$$\overline{P}^{n,\omega}\left(Y_0\neq \overline{Y}_0\right)\leq \left(1/(C_2(\omega)+1)\right)^m+\overline{P}_{n,\omega}\left(\rho^{n,\omega}_m\geq L_n(\omega)\right),$$ 
valid for any positive integer $1\leq m\leq L_n(\omega)-1$. We remind that from Proposition \ref{crucial}, $\lim_{n\rightarrow \infty}L_n(\omega)=\infty$.
Moreover,
$$\overline{P}_{n,\omega}\left(\rho^{n,\omega}_m\geq L_n(\omega)\right)\leq \eta(\omega)^{-L_n(\omega)}\overline{E}_{n,\omega}\left[\eta(\omega)^{\rho^{n,\omega}_m}\right],$$
where $\eta(\omega)>1$ is defined in Lemma \ref{lem3}.  
Writing 
$$\rho^{n,\omega}_m=\rho^{n,\omega}_1+\sum_{j=2}^m\left(\rho^{n,\omega}_j-\rho^{n,\omega}_{j-1}\right)$$
and using Lemma \ref{lem3} and successive conditionings, we get 
$$\overline{P}_{n,\omega}\left(\rho^{n,\omega}_m\geq L_n(\omega)\right)\leq D(\omega)^{m-1}\left(1+V(z)+V(\overline{z})\right)\frac{\eta(\omega)^{m-1}}{\eta(\omega)^{L_n(\omega)}}.$$
We then simply choose, $k=k(\omega)$ such that $\widetilde{\eta}=(D\eta)^{1/k}/\eta<1$ a.s. and set $m=L_n(\omega)/k(\omega)$. It is clear that for a suitable $\widetilde{F}(\omega)>0$ and $\widetilde{\kappa}(\omega)=\max\left(1/C_2(\omega),\widetilde{\eta}(\omega)\right)<1$, we have 
$$\overline{P}^{n,\omega}\left(Y_0\neq \overline{Y}_0\right)\leq \widetilde{F}(\omega)\left(1+V(z)+V(\overline{z})\right)\widetilde{\kappa}(\omega)^{L_n(\omega)}.$$
Since from Proposition \ref{crucial}, $\lim_{n\rightarrow \infty}L_n/n>0$ a.s., we deduce the existence of a pair $(F,\kappa)$ for which the proposed bound follows.$\square$

We now deduce the important result concerning the convergence of the backward iterations of the chains.

\begin{cor}\label{back}
The following assertions hold true.
\begin{enumerate}
\item
The sequence $\left(\delta_zP_{X_{-n}}\cdots P_{X_{-1}}\right)_{n\geq 0}$ is converging $\P-$almost surely in total variation towards a random probability measure denoted by $\pi_{\xi_{-1}}$ and that not depend on $z$. Moreover there exist a pair of random variables $(F,\kappa)$ taking values in $(0,\infty)\times (0,1)$ such that $\P-$a.s.,
\begin{equation}\label{goodconv}
d_{TV}\left(\delta_zP_{X_{-n}}\cdots P_{X_{-1}},\pi_{\xi_{-1}}\right)\leq F\left(1+V(z)\right)\kappa^n.
\end{equation}

\item 
For any $t\in\Z$, we have $\pi_{\xi_{t-1}}P_{X_t}=\pi_{\xi_t}$ a.s. 
\item
If $\left(\nu_t\right)_{t\in\Z}$ is a sequence of identically distributed random probability measures such that $\nu_tP_{X_t}=\nu_{t+1}$ a.s., then 
$\nu_0=\pi_{\xi_0}$ a.s.   
\item
We have $\pi_{\xi_{-1}}V<\infty$ a.s.
\end{enumerate}
\end{cor}

\paragraph{Note.} Let $Y_t^z$, $t\geq 1$ be the successive iterations of the chain (\ref{model}) initialized with $Y_0^z=z$. Let $\pi_t^z$ the probability distribution of $Y_t^z$ and define $\pi^{\infty}(A)=\E\left(\pi_{\xi_{-1}}(A)\right)$ for $A\in\mathcal{B}(E)$. Then by stationarity,
\begin{eqnarray*}\label{compar}
d_{TV}\left(\pi_t^{z},\pi^{\infty}\right)&=&\sup_{A\in \mathcal{B}(E)}\left\vert  \E\left[\delta_zP_{X_0}\cdots P_{X_t}(A)\right]-\pi^{\infty}(A)\right\vert\\
&=& \sup_{A\in\mathcal{B}(E)}\left\vert \E\left[\delta_zP_{X_{-t-1}}\cdots P_{X_{-1}}(A)\right]-\E\left[\pi_{\xi_{-1}}(A)\right]\right\vert\\
&\leq& \E d_{TV}\left(\delta_zP_{X_{-t-1}}\cdots P_{X_{-1}},\pi_{\xi_{-1}}\right).
\end{eqnarray*} 
Using the dominated convergence theorem and the first point of Corollary \ref{back}, we deduce that 
$$\lim_{t\rightarrow \infty}d_{TV}\left(\pi_t^{z},\pi^{\infty}\right)=0.$$
As it will appear more clearly in the proof of Theorem \ref{mainresult}, $\pi^{\infty}$ is simply the marginal distribution of the stationary solution of (\ref{model}).

\paragraph{Proof of Corollary \ref{back}}
We prove this corollary point by point.
\begin{enumerate}
\item
We fix $z\in E$ and for two integers $s<t$, set $Q_{s,t}=P_{X_s}\cdots P_{X_{t-1}}$, we prove that almost surely, $\left(\delta_z Q_{-n,0}\right)_{n\geq 1}$ is a Cauchy sequence on the space of probability measures on $E$, endowed with the total variation metric.
Let $m$ be another integer greater than $n$ and $t_n(\omega)=\tau_{-L_n(\omega)+1}(\omega)$. Observe that $t_n(\omega)-n,t_n(\omega)-m\geq C_1(\omega)$. We have from Assumption {\bf A2}, Proposition \ref{prop1} and Proposition \ref{crucial},
\begin{eqnarray*}
&&d_{TV}\left(\delta_z Q_{-n,0},\delta_z Q_{-m,0}\right)\\
&\leq& \int_E\int_E \delta_zQ_{-n,t_n}(dz')\delta_{\overline{z}'}Q_{-m,t_n}(d\overline{z}')F \left(1+V(z')+V(\overline{z}')\right)\kappa^{-t_n}\\
&\leq& F\left(1+(1-1/C_1)V(z)+2C_1+\prod_{j=-m}^{t_n-1} \lambda\left(X_j\right)V(\overline{z})\right)\kappa^{-t_n}.
\end{eqnarray*}
Observe that we also have $\prod_{j=-m}^{t_n-1} \lambda\left(X_j\right)\leq 1-1/C_1$ but Assumption {\bf A2} guarantees that this product converges also to $0$ as $m\rightarrow \infty$ a.s.
Since, $t_n\rightarrow \infty$ as $n\rightarrow \infty$, we deduce the Cauchy sequence property a.s. and the existence of a limit for $\delta_z Q_{-n,-1}$ that do not depend on $z$, using Proposition \ref{crucial}. Moreover by letting $m\rightarrow \infty$ in the previous bound and using the fact that 
$\lim\inf_{n\rightarrow \infty}(-t_n)/n\geq \lim\inf_{n\rightarrow \infty}(L_n-1)/n>0$ a.s., we deduce the second part of the assertion.

\item
Note first that by stationarity, $\pi_{\xi_t}$ is well defined, has the same probability distribution as $\pi_{\xi_{-1}}$ and is defined by 
$$\pi_{\xi_t}=\lim_{n\rightarrow \infty} \delta_zP_{X_{t-n}}\ldots \P_{X_t},$$
where the limit is almost sure and the convergence in total variation. 
Moreover 
\begin{eqnarray*}
d_{TV}\left(\pi_{X_{t-1}}P_{X_t},\pi_{\xi_t}\right)&\leq& d_{TV}\left(\pi_{\xi_{t-1}} P_{X_t},\delta_z Q_{t-n}^t P_{X_t}\right)+d_{TV}\left(\delta_z Q_{t-n}^t P_{X_t},\pi_{\xi_t}\right)\\
&\leq&d_{TV}\left(\pi_{\xi_{t-1}},\delta_z Q_{t-n}^t\right)+d_{TV}\left(\delta_z Q_{t-n}^t P_{X_t},\pi_{\xi_t}\right)\\
\end{eqnarray*}
and we obtain $d_{TV}\left(\pi_{X_{t-1}}P_{X_t},\pi_{\xi_t}\right)=0$ by letting $n\rightarrow \infty$, which concludes the proof of the second point.
\item
Let $t\in\Z$. For simplicity, we work on the canonical space $\Omega=F^{\Z}$ of the environment and set $X_t(\omega)=\omega_t$ for $\omega\in\Omega$. We also 
denote by $\theta:\Omega\rightarrow \Omega$ the translation operator defined by $\theta\cdot\omega=(\omega_{t+1})_{t\in\Z}$. 
It is sufficient to show that for any $\epsilon>0,\eta>0$, $\P\left[d_{TV}\left(\nu_{t-1},\pi_{\xi_{t-1}}\right)>\eta\right]\leq \epsilon$.
Since for $n\geq 1$, $\nu_{t-1}=\nu_{t-n-1}P_{X_{t-n}}\cdots P_{t-1}$, we get 
\begin{eqnarray*}
&&d_{TV}\left(\nu_{t-1},\pi_{\xi_{t-1}}\right)\\
&\leq& d_{TV}\left(\nu_{t-n-1}P_{X_{t-n}}\cdots P_{X_{t-1}},\delta_z P_{X_{t-n}}\cdots P_{X_{t-1}}\right)+d_{TV}\left(\delta_z P_{X_{t-n}}\cdots P_{X_{t-1}},\pi_{\xi_{t-1}}\right)\\
&=& A_n+B_n.
\end{eqnarray*}
We have, using Proposition \ref{prop1},
\begin{eqnarray*} 
A_n&\leq& \int \nu_{t-n-1}(d\overline{z})d_{TV}\left(\delta_{\overline{z}} P_{X_{t-n}}\cdots P_{X_{t-1}},\delta_z P_{X_{t-n}}\cdots P_{X_{t-1}}\right)\\
&\leq& \nu_{t-n-1}\left(\{V>R\}\right)+\left(1+V(z)+R\right)\cdot F\circ \theta^t\cdot \left(\kappa\circ \theta^t\right)^n.
\end{eqnarray*}
Take $R>0$ such that 
$$\P\left(\nu_{t-n-1}\left(\{V>R\}\right)>\eta/3\right)=\P\left(\nu_0\left(\{V>R\}\right)>\eta/3\right)\leq \epsilon/3.$$
This is possible since from the finiteness of $V$, $\lim_{R\rightarrow \infty}\nu_0\left(\{V>R\}\right)=0$ a.s.
For such $R$, choose now $n$ sufficiently large in order to have
$$\P\left(\left(1+V(z)+R\right)\cdot F\circ \theta^t\cdot \left(\kappa\circ \theta^t\right)^n>\eta/3\right)\leq \epsilon/3$$
and also $\P(B_n>\eta/3)\leq \epsilon/3$. This leads to the result.
\item
For $M>0$, let $V\wedge M:E\rightarrow \R$ be the function defined by $(V\wedge M)(y)=\min\left(V(y),M\right)$ for $y\in E$.
We have for $z\in E$ and any $M>0$, 
\begin{eqnarray*}
\delta_z P_{X_{-n}}\cdots P_{X_{-1}}(V\wedge M)&\leq& \delta_z P_{X_{-n}}\cdots P_{X_{-1}}V\\
&\leq& \prod_{i=1}^n\lambda\left(X_{-i}\right)V+b\left(X_{-1}\right)+\sum_{i\geq 2}\prod_{j=1}^{i-1}\lambda\left(X_{-j}\right)b\left(X_{-i}\right)
\end{eqnarray*}
and using {\bf A2}, we have already shown that the last bound is finite a.s. and that 
$$\lim_{n\rightarrow \infty}\prod_{i=1}^n\lambda\left(X_{-i}\right)=0\mbox{ a.s.}.$$
By letting $n\rightarrow \infty$ and using (\ref{goodconv}), we get $\pi_{\xi_{-1}}V=\sup_{M>0}\pi_{\xi_{-1}}(V\vee M)<\infty$ a.s. which leads to the result.
\end{enumerate}  

\subsection{Proof of Theorem \ref{mainresult}}
We proceed as in \citet{DNT}, using Proposition \ref{prop1}. As in the proof of the third point of Proposition \ref{prop1}, we consider the canonical space $\Omega=F^{Z}$ of the environment and we set $\nu^{\omega}=\pi_{\xi_{-1}(\omega)}$.
For proving existence of a stationary process, we consider the finite-dimensional distributions 
\begin{equation}\label{finite}
\zeta_{u,t}^{\omega}(dy_u,\cdots,dy_t)=\nu^{\theta^u\omega}\left(dy_u\right)\prod_{i=u}^{t-1}P_{\omega_i}\left(y_i,dy_{i+1}\right),
\end{equation}
for $u\leq t$ in $\Z$. Using Kolmogorov's extension theorem, there exists a unique probability measure $\zeta^{\omega}$ on $E^{\Z}$ compatible with such a family.
On $\widetilde{\Omega}=E^{\Z}\times \Omega$, the probability measure $d\widetilde{P}\left(y,\omega\right)=d\zeta^{\omega}(y)d\P(\omega)$ is solution of (\ref{model}), in the sense that if $Y_t(y,\omega)=y_t$ and $X_t(y,\omega)=\omega_t$ for $t\in\Z$, we have
$$\widetilde{\P}\left(Y_t\in A\vert X,Y_{t-1},Y_{t-2},\ldots\right)=P_{X_{t-1}}\left(Y_{t-1},A\right).$$
Stationarity of the process $\left((X_t,Y_t)\right)_{t\in \Z}$ follows from the equalities
$$\zeta_{u,t}^{\theta^n \omega}=\zeta_{u+n,t+n}^{\omega},$$
for $u\leq t$ in $\Z$ and $n\in \N$. Let us now show uniqueness. If a stationary stochastic process $\left((Y_t,X_t)\right)_{t\in\Z}$ defined on $\left(E\times F\right)^{\Z}$ and with probability distribution $\mathbb{Q}$  satisfies (\ref{model}), we will have almost surely,
$$\mathbb{Q}\left(Y_t\in A\vert X\right)=\nu^{\theta^t X}(A),\quad A\in \mathcal{B}(E).$$
Indeed, if $\overline{\nu}^{X}_t(A)=\mathbb{Q}\left(Y_t\in A\vert X\right)$, then $\overline{\nu}^{X}_t$ is a random probability measure such that $\overline{\nu}^{X}_tP_{X_{t+1}}=\overline{\nu}^{X}_{t+1}$ $\P-$a.s.
Moreover, by stationarity of the process, we have $\overline{\nu}^{X}_t=\overline{\nu}^{\theta^t X}_0$ a.s. and $\overline{\nu}^{X}_t$ has then the same probability distribution has $\overline{\nu}^X_{t+1}$. Using the third point of Proposition \ref{prop1}, we get 
$\overline{\nu}^{X}_t=\nu^{\theta^t X}$ a.s. The finite-dimensional distributions of $\mathbb{Q}$ and $\widetilde{\P}$ are then automatically the same and we conclude that $\mathbb{Q}=\widetilde{\P}$.

We now prove ergodicity of the unique stationary solution.
As in \citet{DNT}, one can use the uniqueness property derived in the third point of Proposition \ref{prop1} to get the ergodicity property.
However, we will provide another one, which is more direct and analog to the proof given by \citet{Kifer}.
 
Denoting by $\theta_{*}$ the translation operator on $E^{N}$, i.e. $\theta_{*}y=(y_{t+1})_{t\in\N}$ for every $y\in E^{\N}$, one can define the application 
$\tau=\left(\theta_{*},\theta\right)$ on $\widetilde{\Omega}$ by $\tau((y,\omega))=\left(\theta_{*}y,\theta \omega\right)$ for $(y,\omega)\in \widetilde{\Gamma}$.
Let $p:E^{\Z}\rightarrow E^{\N}$ the projection operator defined by $p\left((y_t)_{t\in \Z}\right)=(y_t)_{t\in\N}$.
If $\omega\in \Omega$, we define a probability measure $\gamma_{*}^{\omega}$ on $E^{\N}$ by $\gamma_{*}^{\omega}(A)=\gamma^{\omega}\left(p^{-1}(A)\right)$ for any Borel subset $A$ of $E^{\N}$. Finally, we define a probability measure $\gamma_{*}$ on 
$\widetilde{\Omega}_{*}:=E^{\N}\times \Omega$ by $\gamma_{*} f=\int f(y,\omega) \gamma_{*}^{\omega}(dy)\P(d\omega)$ for every measurable and bounded function $f:\widetilde{\Omega}_{*}\rightarrow \R$. We are going to show that $\tau$ is ergodic for $\gamma_{*}$, i.e. if $I\subset \widetilde{\Omega}_{*}$ is a measurable invariant set for $\tau$, i.e. $\tau^{-1}I=I$, then $\gamma_{*}(I)\in \{0,1\}$. Note that $\gamma_{*}$ is simply the probability distribution 
of the pair $\left((Y_t)_{t\in \N},(X_t)_{t\in\Z}\right)$ under $\widetilde{\P}$ where $\left(Y_t,X_t)\right)_{t\in\Z}$ is the stochastic process defined just above. Let us first explain why this ergodicity property entails ergodicity of the process $\left(Z_t\right)_{t\in\N}$ with $Z_t=(Y_t,X_t)$. 
To this end, let $J\subset (E\times F)^{\N}$ an invariant subset for the translation operator $\theta_{**}$ on $(E\times F)^{\N}$. The process $(Z_t)_{t\in\N}$ 
is ergodic if and only if $\widetilde{\P}\left((Z_t)_{t\in\N}\in J\right)\in \{0,1\}$ for every measurable invariant set, i.e. $\theta_{**}^{-1}J=J$.
But from the definition of $\widetilde{\P}$, we have 
$$\widetilde{\P}\left(Z\in J\right)=\gamma_{*}(I),$$
where 
$$I=\left\{(y,\omega)\in E^{\N}\times \Omega: ((y_t,\omega_t))_{t\in\N}\in J\right\}.$$
Moreover
$$(y,\omega)\in I\Leftrightarrow \tau\cdot(y,\omega)\in I,$$
meaning that $I$ is an invariant set for $\tau$. If we show that $\tau$ is ergodic for $\gamma_{*}$, we will then get $\widetilde{P}\left(Z\in J\right)\in \{0,1\}$ and then $\theta_{**}$ will be ergodic for $\widetilde{P}_{Z}$ which is precisely the definition of the ergodicity of the process $(Z_t)_{t\in\N}$.
From this one-sided ergodicity property, one can also deduce the two-sided ergodicity property for the process $(Z_t)_{t\in\Z}$, which means that 
$\widetilde{\P}\left((Z_t)_{t\in\Z}\in K\right)\in \{0,1\}$ for every invariant subset $K$ of $(E\times F)^{\Z}$, that is $K$ is a measurable subset of $(E\times F)^{\Z}$ such that $\theta^{-1}K=K$ and $\theta\cdot z=(z_{t+1})_{t\in\Z}$ for every $z\in (E\times F)^{\Z}$. 
A proof explaining why the one-sided ergodicity implies the two-sided ergodicity property can be found in \citet{DDM}, Theorem $33$.

Now, let $I$ be an invariant set for $\tau$, i.e. $\theta_{*}^{-1} I=I$. For $\omega\in \Omega$, let $I^{\omega}$ be the section of $\Gamma$, i.e.
$$I^{\omega}=\left\{y\in E^{\N}: (y,\omega)\in I\right\}.$$
From the invariance property of $\Gamma$, we have the equalities 
\begin{equation}\label{invariance}
I^{\omega}=\theta_{1,*}^{-n}I^{\theta^n\omega},\quad n\in\N^{*}.
\end{equation}
For any measurable subset $A$ included in $E^{\N}$, we get from the definition of $\gamma^{\omega}$ the equalities
$$\gamma_{*}^{\theta^n\omega}(A)=\gamma^{\omega}_{*}\left(\theta_{*}^{-n}A\right),\quad n\in\N^{*}.$$
Using (\ref{invariance}), we deduce that the equalities
\begin{equation}\label{invariance2}
\gamma^{\omega}_{*}\left(I^{\omega}\right)=\gamma_{*}^{\omega}\left(\theta_{1,*}^{-n}I^{\theta^n\omega}\right)=\gamma_{*}^{\theta^n \omega}\left(I^{\theta^n\omega}\right),\quad n\in\N^{*}.
\end{equation}
We are going to show that there exist some random variables $n_i:\Omega\rightarrow \N^{*}$ such that $n_i<n_{i+1}$ for every $i\in \N^{*}$ and 
\begin{equation}\label{invariance3}
\lim_{i\rightarrow \infty}\gamma^{\omega}_{*}\left(I^{\omega}\cap \theta_{*}^{-n_i(\omega)}I^{\theta^{n_i(\omega)}\omega}\right)=\gamma_{*}^{\omega}\left(I^{\omega}\right)^2.
\end{equation}
Note that (\ref{invariance3})and (\ref{invariance2}) automatically entails that $f(\omega):=\gamma_{*}^{\omega}\left(I^{\omega}\right)\in\{0,1\}$ and since $f(\theta\omega)=f(\omega)$, ergodicity of $\theta$ for the measure $\P$ implies that $f$ is $\P-$almost surely constant and then that 
$$\gamma_{*}(I)=\int_{\Omega}\gamma_{*}^{\omega}\left(I^{\omega}\right)\P(d\omega)\in\{0,1\},$$
which will prove the ergodicity of $\tau$ for $\gamma_{*}$.
To prove (\ref{invariance3}), we will use (\ref{goodconv}). Let $M>1$ such that $\P\left(F\leq M, \kappa\leq 1-1/M\right)>0$. 
This is always possible since $(F,\kappa)\in (0,\infty)\times (0,1)$ $\P-$a.s.
We then define the $n_i$'s as the successive return times of the pair $(F\circ \theta^n,\kappa\circ\theta^n)$ in the Cartesian product $[0,M]\times [0,1/1/M]$.
From ergodicity of $\theta$ for $\P$, these return times are always finite. 
Next we consider $A$ and $B$ two Borel subsets of $E^{\N}$. We assume that $A$ is a cylinder set, i.e. $A=\prod_{j\in\N}A_j$ with $A_j\in \mathcal{B}(E)$ and $A_j=E$ for $j\geq k+1$. 
For $y_n\in E$, let also $B_{y_n}^{+}=\left\{z\in E^{\N}: (y_n,z)\in B\right\}$. For simplicity of notations, we set for $s<t$,$t$ being possibly infinite,
$$\nu^{\omega}_{s,t}\left(y_s,d(y_{s+1},\ldots,y_t)\right)=\otimes_{i=s}^{t-1}P_{\omega_i}(y_i,dy_{i+1})$$
and when $t$ is finite, $Q^{\omega}_{s,t}(y_s,dy_t)=\delta_{y_s}\prod_{i=s}^{t-1}P_{\omega_i}(dy_t)$.
We have for a $n\in\N^{*}$,
\begin{eqnarray*}
&&\gamma_{*}^{\omega}\left(A\cap \theta_{*}^{-n} B\right)\\
&=& \int_{A_0\times\cdots A_k}\pi^{\xi_{-1}(\omega)}(dy_0)\nu^{\omega}_{0,k}\left(y_0,d(y_1,\ldots,y_k)\right)\delta_{y_k}Q^{\omega}_{k,n}(dy_n)\nu_{n,\infty}^{\omega}\left(y_n,B_{y_n}^{+}\right),
\end{eqnarray*}
On the other hand, we have for $n>k$,
\begin{eqnarray*}
&&\gamma_{*}^{\omega}(A)\gamma_{*}^{\omega}\left(\theta_{*}^{-n}B\right)\\
&=& \gamma_{*}^{\omega}(A)\gamma_{*}^{\theta^n\omega}(B)\\
&=& \int_{A_0\times\cdots A_k}\pi^{\xi_{-1}(\omega)}(dy_0)\nu^{\omega}_{0,k}\left(y_0,d(y_1,\ldots,y_k)\right)\pi_{\xi_{n-1}(\omega)}(dy_n)\nu_{n,\infty}^{\omega}\left(y_n,B_{y_n}^{+}\right)
\end{eqnarray*}
and using (\ref{goodconv}), we get 
\begin{eqnarray*}
&& \left\vert \gamma^{\omega}_{*}\left(A\cap \theta_{*}^{-n}B\right)-\gamma^{\omega}_{*}(A)\gamma^{\omega}_{*}\left(\theta_{*}^{-n}B\right)\right\vert\\
&\leq& \int_E \pi_{\xi_{k-1}(\omega)}(dy_k)\left(1+V(y_k)\right)F\left(\theta^n\omega\right)\kappa\left(\theta^n\omega\right)^{n-k}.
\end{eqnarray*}
Using the definition of the random times $n_i$ and Corollary \ref{back}(4.), we deduce that for $\P-$almost every $\omega$, we have for all cylinder set $A$,
\begin{equation}\label{goodinter}
\lim_{i\rightarrow \infty}\sup_{B\in \mathcal{B}(E^{\N})}\left\vert \gamma_{*}^{\omega}\left(A\cap \theta_{*}^{-n_i(\omega)}B\right)-\gamma_{*}^{\omega}(A)\gamma_{*}^{\omega}\left(\theta_{*}^{-n_i(\omega)}B\right)\right\vert=0.
\end{equation}
By approximating $A$ by a finite union of disjoint cylinder sets for the measure $\gamma_{*}^{\omega}$, one can get (\ref{goodinter}) for any Borel subset $A$ of $E^{\N}$. Applying (\ref{goodinter}) to $A=I^{\omega}$, we deduce (\ref{invariance3}) and then the ergodicity of $\tau$ for $\gamma$. As explained above, this lead to the ergodicity property of Theorem \ref{mainresult}.$\square$

\subsection{Proof of Theorem \ref{mainresult2}}

We apply the machinery for the proof of Theorem \ref{mainresult} using a blocking argument. We only discuss the required modifications.
For $j=0,\ldots,p-1$ and $i\in\Z$, let 
$$X_i^{(j)}=\left(X_{j+(i-1)p+1},\ldots,X_{j+ip}\right).$$
We replace the environment $(X_t)_{t\in\Z}$ by $(X_t^{(j)})_{t\in\Z}$ and set 
$$P_{X_t^{(j)}}=\prod_{s=j+(t-1)p+1}^{j+tp}P_{X_s}.$$
From our assumptions, Corollary \ref{back} applies to this family of random kernels and setting $\xi_i^{(j)}=\left(X^{(j)}_s\right)_{s\leq i}$, we deduce the existence of a random probability measure $\pi_{\xi^{(j)}-1}^{(j)}$ such that for any $z\in E$,
$$d_{TV}\left(\delta_zP_{X^{(j)}_{-n}}\cdots P_{X^{(j)}_{-1}},\pi^{(j)}_{\xi^{(j)}_{-1}}\right)\leq F_j\left(1+V(z)\right)\kappa_j^n,$$
where $F_j:\Omega\rightarrow (0,\infty)$ and $\kappa_j:\Omega\rightarrow (0,1)$ are random variables. 
Proceeding as in the proof of Corollary \ref{back}, point $3.$, it is not difficult to show that $\pi^{(j)}_{\xi_{-1}^j}P_{X_{-p+j+1}}=\pi^{(j+1)}_{\xi_{-1}^{j+1}}$ a.s. for $j=0,\ldots, p-2$.
Finally, setting $\pi_{\xi_{-1}}=\pi^{(p-1)}_{\xi^{(p-1)}_{-1}}$, one can deduce that (\ref{goodconv}) holds true, as well as all the other statements of Corollary \ref{back}. 
The rest of the proof is identical to that of Theorem \ref{mainresult}.

\subsection{Proof of Proposition \ref{ex1}}
Setting for $(x,y)\in \R^d\times \R$,
$$s(x,y)=\left(b_1(x)+a_1(x)y\right)\mathds{1}_{y\leq r(x)}+\left(b_2(x)+a_2(x)y\right)\mathds{1}_{y>r(x)},$$
we apply Theorem \ref{mainresult} setting $P_x(y,A)=\int_A f(y'-s(x,y))dy$ and $V(y)=\vert y\vert$. 
It is straightforward to get {\bf A2} with constants $\lambda(x)$ and $b(x)=\vert b_1(x)\vert+\vert b_2(x)\vert+\E\vert\varepsilon_0\vert$.
If $R>0$, Assumption {\bf A3} is satisfied with $\nu_R(x,A)=(2R)^{-1}\lambda(A\cap [-R,R])$ and $\eta(x,R)=\inf_{\vert y\vert,\vert y'\vert\leq R}f\left(y'-s(x,y)\right)$ which is positive since if $\vert y\vert,\vert y'\vert\leq R$,
$$\vert y'-s(x,y)\vert\leq J(x):=R+\max_{i=1,2}\vert b_i(x)\vert+\max_{i=1,2}\vert a_i(x)\vert R$$
and one can use the lower bound of $f$ of $[-J(x),J(x)]$ to conclude.$\square$

\subsection{Proof of Proposition \ref{ex2}}

We only check the assumptions of Theorem \ref{mainresult2} for the case $p=2$. The general case follows from similar arguments but using more tedious notations.
We set $g(x,y_1,y_2)=\sum_{j=1}^1a_j(x,y_1,y_2)y_j$ and $V(y_1,y_2)=\vert y_1\vert+\vert y_2\vert$. We denote by $\Vert\cdot\Vert$ the $\ell_1-$norm on $\R^2$
and we denote in the same way the corresponding operator norm for two-dimensional square matrices.
First, we note that it is equivalent to study the following two-dimensional autoregressive process of order one:
$$\begin{pmatrix} Y_t\\ Y_{t-1}\end{pmatrix}=\begin{pmatrix} g\left(X_{t-1},Y_{t-1},Y_{t-2}\right)+\varepsilon_t\\Y_{t-1}\end{pmatrix}.$$
Note also that for this bivariate dynamic, we have
$$P_x((y_1,y_2),A)\int \mathds{1}_A\left(g(x,y_2,y_1)+u,y_2\right)f(u)du,\quad (y_1,y_2,x,A)\in \R^{d+1}\times \mathcal{B}(\R).$$
Setting $Z_t=(Y_t,Y_{t-1})'$, $H(x,y_1,y_2)=(g(x,y_1,y_2),0)'$ and $B_t=(\varepsilon_t,0)'$, we then have
$$\vert Z_t\vert \leq A\left(X_{t-1}\right)\cdot\vert Z_{t-1}\vert+\vert B_t\vert,$$
where $\vert z \vert=(\vert z_1\vert, \vert z_2\vert)'$ and $\leq $ stands for the componentwise ordering.
Iterating this bound, we have 
$$\vert Z_t\vert\leq A(X_{t-1})\cdots A(X_{t-k})\vert Z_{t-k}\vert+\vert B_t\vert+\sum_{j=2}^k\prod_{i=1}^{j-1}A(X_{t-i})\vert B_{t-j-1}\vert.$$
We then get the bound
$$P_{x_1}\cdots P_{x_k}V\leq \lambda(x_k,\ldots, x_1)V+b(x_k,\ldots,x_1),$$
with $\lambda(x_k,\ldots, x_1)=\Vert A(x_k)\cdots A(x_1)\Vert$ and 
$$b(x_k,\ldots,x_1)=\Vert c\Vert+\sum_{j=2}^k\Vert A(x_k)\Vert\cdots\Vert A(x_{k-j+2})\Vert\cdot\Vert c\Vert$$
with $c=(\E\vert \varepsilon_0\vert,0)'$.
Assumption {\bf A4} is then satisfied if $k$ is large enough, since under the mixing assumption, the limit in (\ref{multi}) is given by
$\exp\left(\E\log\Vert A(X_{-1})\cdots A(X_{-k}\Vert\right)$ which is less than one from the negativity of the Lyapunov exponent.
One can choose $p=k=2k'$ with $k'$ large enough. 

It remains to check {\bf A5}.  
To this end, if $V(y_1,y_2)\leq R$ and $A\in\mathcal{B}(\R^2)$, we define $\nu_R(\cdot)$ as the Lebesgue measure over $[-R,R]^2$. We have 
\begin{eqnarray*}
\left[P_{x_1}P_{x_2}\right]((y_1,y_2),A)&=&\int_A f\left(v_1-g(x_1,y_2,y_1)\right)f\left(v_2-g(x_2,v_1,y_2)\right)dv_1dv_2\\
&\geq & \widetilde{\eta}(R,x_2,x_1)\nu_R(A), 
\end{eqnarray*}
with 
$$\widetilde{\eta}(R,x_2,x_1)=(2R)\inf_{\vert y_1\vert,\vert y_2\vert,\vert v_1\vert,\vert v_2\vert\leq R}f\left(v_1-g(x_1,y_2,y_1)\right)f\left(v_2-g(x_2,v_1,y_2)\right)>0.$$
We then get
$$\left[\prod_{i=1}^p P_{x_i}\right]((y_1,y_2),A)\geq \nu_R\left(\{V\leq R\}\right)^{k'-1}\prod_{i=1}^{k'}\widetilde{\eta}(R,x_{2i},x_{2i-1})\nu_R(A),$$
which shows that $\{V\leq R\}$ is always a small set in the sense of {\bf A5}.
The result then follows from Theorem \ref{mainresult2}.$\square$

\section{Appendix: proof of Lemma \ref{lem3}}\label{Section6}
For simplicity of notations and since $\omega\in\Omega$ is fixed, we omit the subscripts $\omega$ and $n$ in all our notations.
Set $x=V(z)+V(\overline{z})$ and assume that $x>R$. Let $\mathcal{H}_i=\sigma\left((Z_j,\overline{Z}_j): j\leq i\right)$.
Remind the third inequality given in Lemma \ref{lem2}. Since $(1-1/C_1)y+2C_1\leq \eta^{-1}y-1$ when $y>R$, we get 
$$E\left(W_{i+1}\vert \mathcal{H}_i\right)\leq (1-1/C_1)V(Z_i)+2C_1\leq \eta^{-1}V(Z_i)-1$$
on the event $\{V(Z_{i-1})>R\}$.
This yields to 
\begin{equation}
W_i-\eta \E\left(W_{i+1}\vert \mathcal{H}_i\right)\geq \eta\mbox{ on the event } \{W_i>R\}.
\end{equation}
\begin{enumerate}
\item
We first prove the first point.
We have for $N\geq 1$,
\begin{eqnarray*}
&&\sum_{i=0}^N \eta^{i+1}(\rho_1\geq i+1)\\
&\leq& \sum_{i=0}^N E\left[\mathds{1}_{\rho_1\geq i+1}\left(\eta^i W_i-\eta^{i+1}E\left(W_{i+1}\vert \mathcal{H}_i\right)\right)\right]\\
&=&\sum_{i=0}^N E\left[\mathds{1}_{\rho_1\geq i+1}\left(\eta^i W_i-\eta^{i+1}W_{i+1}\right)\right]\\
&\leq& \sum_{i=0}^N E\left[\mathds{1}_{\rho_1\geq i}\eta^i W_i-\mathds{1}_{\rho_1\geq i+1}\eta^{i+1}W_{i+1}\right]\\
&=& x-\eta^{N+1}E\left[\mathds{1}_{\rho_1\geq N+1}W_{N+1}\right]\leq x.
\end{eqnarray*}
The first inequality follows from the fact that on the event $\{\rho_1\geq i+1\}$, we have $W_i>R$.
We then deduce that $\sum_{i\geq 0}\eta^{i+1}P\left(\rho_1\geq i+1\right)$ is finite then that $\rho_1$ is finite a.s. Moreover we have 
$$E\left(\eta^{\rho_1}\right)\leq \sum_{i\geq 0}\eta^{i+1}P\left(\rho_1>i+1\right)\leq x,$$
which concludes the proof of the first point
\item
Suppose now that $\rho_m<\infty$ a.s. for some positive integer $m$. Let $A\in \mathcal{H}_{\rho_m}$, $k$ a positive integer and set 
$$B_{k+1}=\{W_{k+1}>R\}\cap A\cap\{\rho_m=k\}.$$
Observe that $B_{k+1}$ is an event of $\mathcal{H}_{k+1}$. Setting $r_{m+1}=\rho_{m+1}-\rho_m$, we first note that on the event $\{r_{m+1}\geq \ell+1,\rho_m=k\}$,
we have $W_{k+\ell}>R$. Then, we get for an integer $N\geq 1$,
\begin{eqnarray*}
&&\sum_{\ell=1}^N\eta^{\ell+1}P\left(\{r_{m+1}\geq \ell+1\}\cap B_{k+1}\right)\\
&\leq& \sum_{\ell=1}^NE\left[\left(\eta^{\ell}W_{k+\ell}-\eta^{\ell+1}E\left(W_{k+\ell+1}\vert \mathcal{H}_{k+\ell}\right)\right)\mathds{1}_{r_{m+1}\geq \ell+1}\mathds{1}_{B_{k+1}}\right]\\
&\leq& \sum_{\ell=1}^NE\left[\left(\eta^{\ell}W_{k+\ell}-\eta^{\ell+1}W_{k+\ell+1}\right)\mathds{1}_{r_{m+1}\geq \ell+1}\mathds{1}_{B_{k+1}}\right]\\
&\leq &\sum_{\ell=1}^NE\left[\left(\eta^{\ell}W_{k+\ell}\mathds{1}_{r_{m+1}\geq \ell}-\eta^{\ell+1}W_{k+\ell+1}\mathds{1}_{r_{m+1}\geq \ell+1}\right)
\mathds{1}_{B_{k+1}}\right]\\
&=& E\left[\left(\eta W_{k+1}-\eta^{N+1}W_{k+N+1}\mathds{1}_{r_{m+1}\geq N+1}\right)\mathds{1}_{B_{k+1}}\right]\\
&\leq& \E\left[\eta W_{k+1}\mathds{1}_{B_{k+1}}\right]\\
&\leq & \E\left[\eta \E\left[W_{k+1}\vert \mathcal{H}_k\right]\mathds{1}_{A\cap  \{\rho_m=k\}}\right]\\
&\leq &\eta \E\left[\left((1-1/C_1)V(W_k)+2 C_1\right)\mathds{1}_{A\cap  \{\rho_m=k\}}\right]\\
&=& \eta\left((1-1/C_1)R+2C_1\right)\P\left(A\cap \{\rho_m=k\}\right).
\end{eqnarray*}
Moreover, setting $B_{k+1}'=A\cap\{\rho_m=k\}\cap \{W_{k+1}\leq R\}$, we have 
$$\sum_{\ell=0}^{\infty}\eta^{\ell+1}P\left(\{r_{m+1}\geq \ell+1\}\cap B'_{k+1}\right)=\eta P\left(B'_{k+1}\right)\leq \eta P\left(A\cap\{\rho_m=k\}\right).$$
Our computations lead to the bound
$$\sum_{\ell\geq 0}\eta^{\ell+1}P\left(\{r_{m+1}\geq \ell+1\}\cap A\cap \{\rho_m=k\}\right)\leq \eta \left(1+(1-1/C_1)R+2C_1\right)P\left(A\cap \{\rho_m=k\}\right).$$
Summing over $k$ and taking $A=E$, we get
$$E\left(\eta^{r_m}\right)\leq \sum_{\ell\geq 0}\eta^{\ell+1}P\left(r_{m+1}\geq \ell+1\right)<\infty$$
and then $r_m<\infty$ a.s. This also yields to
$$E\left(\eta^{r_m}\mathds{1}_{A\cap\{\rho_m=k\}}\right)\leq \left(1+(1-1/C_1)R+2C_1\right)P\left(A\cap \{\rho_m=k\}\right)$$
Summing over $k$, we get 
$$\E\left(\eta^{r_m}\mathds{1}_A\right)\leq \left(1+(1-1/C_1)R+2C_1\right)P\left(A\right)$$
and since $A\in \mathcal{H}_{\rho_m}$ is arbitrary, we get the desired result.$\square$
\end{enumerate}
\bibliographystyle{plainnat}
\bibliography{biblogistic}

\begin{thebibliography}{18}
\providecommand{\natexlab}[1]{#1}
\providecommand{\url}[1]{\texttt{#1}}
\expandafter\ifx\csname urlstyle\endcsname\relax
  \providecommand{\doi}[1]{doi: #1}\else
  \providecommand{\doi}{doi: \begingroup \urlstyle{rm}\Url}\fi

\bibitem[Bougerol and Picard(1992)]{Boug}
P.~Bougerol and N.~Picard.
\newblock {Stationarity of GARCH processes and of some nonnegative time
  series}.
\newblock \emph{Journal of econometrics}, 52\penalty0 (1-2):\penalty0 115--127,
  1992.

\bibitem[Cai et~al.(2000)Cai, Fan, and Yao]{Fan}
Zongwu Cai, Jianqing Fan, and Qiwei Yao.
\newblock Functional-coefficient regression models for nonlinear time series.
\newblock \emph{Journal of the American Statistical Association}, 95\penalty0
  (451):\penalty0 941--956, 2000.

\bibitem[Chamberlain(1982)]{Chamb}
G.~Chamberlain.
\newblock {The general equivalence of Granger and Sims causality}.
\newblock \emph{Econometrica: Journal of the Econometric Society}, pages
  569--581, 1982.

\bibitem[Chen and Tsay(1993)]{chen}
Rong Chen and Ruey~S Tsay.
\newblock Functional-coefficient autoregressive models.
\newblock \emph{Journal of the American Statistical Association}, 88\penalty0
  (421):\penalty0 298--308, 1993.

\bibitem[Cogburn(1984)]{Cogburn}
R.~Cogburn.
\newblock {The ergodic theory of Markov chains in random environments}.
\newblock \emph{Zeitschrift f{\"u}r Wahrscheinlichkeitstheorie und verwandte
  Gebiete}, 66\penalty0 (1):\penalty0 109--128, 1984.

\bibitem[Douc et~al.(2013)Douc, Doukhan, and Moulines]{DDM}
R.~Douc, P.~Doukhan, and E.~Moulines.
\newblock Ergodicity of observation-driven time series models and consistency
  of the maximum-likelihood estimator.
\newblock \emph{Stochastic Processes and their Applications}, 123:\penalty0
  2620--2647, 2013.

\bibitem[Douc et~al.(2004)Douc, Moulines, and Rosenthal]{DM}
Randal Douc, Eric Moulines, and Jeffrey~S Rosenthal.
\newblock Quantitative bounds on convergence of time-inhomogeneous markov
  chains.
\newblock \emph{Annals of Applied Probability}, pages 1643--1665, 2004.

\bibitem[Douc et~al.(2014)Douc, Moulines, and Stoffer]{DoucMoulines}
Randal Douc, Eric Moulines, and David Stoffer.
\newblock \emph{Nonlinear time series: Theory, methods and applications with R
  examples}.
\newblock CRC press, 2014.

\bibitem[Doukhan and Neumann(2019)]{DN}
P.~Doukhan and M.~H. Neumann.
\newblock Absolute regularity of semi-contractive garch-type processes.
\newblock \emph{Journal of Applied Probability}, 56\penalty0 (1):\penalty0
  91--115, 2019.

\bibitem[Doukhan et~al.(2020)Doukhan, Neumann, and Truquet]{DNT}
Paul Doukhan, Michael~H Neumann, and Lionel Truquet.
\newblock Stationarity and ergodic properties for some observation-driven
  models in random environments.
\newblock \emph{arXiv preprint arXiv:2007.07623}, 2020.

\bibitem[Kifer(1996)]{Kifer}
Y.~Kifer.
\newblock {Perron-Frobenius theorem, large deviations, and random perturbations
  in random environments}.
\newblock \emph{Mathematische Zeitschrift}, 222\penalty0 (4):\penalty0
  677--698, 1996.

\bibitem[Lovas and R{\'a}sonyi(2021)]{Lovas}
Attila Lovas and Mikl{\'o}s R{\'a}sonyi.
\newblock Markov chains in random environment with applications in queuing
  theory and machine learning.
\newblock \emph{Stochastic Processes and their Applications}, 137:\penalty0
  294--326, 2021.

\bibitem[Orey(1991)]{Orey}
S.~Orey.
\newblock Markov chains with stochastically stationary transition
  probabilities.
\newblock \emph{The Annals of Probability}, 19\penalty0 (3):\penalty0 907--928,
  1991.

\bibitem[Rosenthal(1995)]{Ro}
Jeffrey~S Rosenthal.
\newblock Minorization conditions and convergence rates for markov chain monte
  carlo.
\newblock \emph{Journal of the American Statistical Association}, 90\penalty0
  (430):\penalty0 558--566, 1995.

\bibitem[Samorodnitsky(2016)]{sam}
Gennady Samorodnitsky.
\newblock \emph{Stochastic processes and long range dependence}, volume~26.
\newblock Springer, 2016.

\bibitem[Stenflo(2001)]{Stenflo}
{\"O}.~Stenflo.
\newblock Markov chains in random environments and random iterated function
  systems.
\newblock \emph{Transactions of the American Mathematical Society},
  353\penalty0 (9):\penalty0 3547--3562, 2001.

\bibitem[Tong(2012)]{tong}
H.~Tong.
\newblock \emph{Threshold models in non-linear time series analysis},
  volume~21.
\newblock Springer Science \& Business Media, 2012.

\bibitem[Tsay(1989)]{Tsay}
R.S. Tsay.
\newblock Testing and modeling threshold autoregressive processes.
\newblock \emph{Journal of the American statistical association}, 84\penalty0
  (405):\penalty0 231--240, 1989.

\end{thebibliography}

\end{document}